\documentclass[12pt]{amsart}
\usepackage{amsfonts,latexsym,rawfonts,amsmath,amssymb,amsthm,mathrsfs,esint}
\usepackage[plainpages=false]{hyperref}
\usepackage{graphicx}
\usepackage{comment}
\usepackage{pgfplots}

\numberwithin{equation}{section}

\RequirePackage{color}
\theoremstyle{plain}

\newtheorem{remark}{Remark}[section]
\newtheorem{theorem}{Theorem}[section]
\newcommand{\beq}{\begin{equation}}
\newcommand{\eeq}{\end{equation}}
\newcommand{\beqs}{\begin{eqnarray*}}
\newcommand{\eeqs}{\end{eqnarray*}}
\newcommand{\beqn}{\begin{eqnarray}}
\newcommand{\eeqn}{\end{eqnarray}}
\newcommand{\beqa}{\begin{array}}
\newcommand{\eeqa}{\end{array}}

\def\S{\mathbb S}
\def\R{\mathbb R}

\def\K{\mathcal K}

\def\vol{\text{Vol}}

\def\s{\sigma_{\S^n}}

\def\Om{\Omega}

\begin{document}
%\linenumbers
\title{On a generalised Blaschke-Santal\`{o} Inequality}

\author{Haodi Chen}
\address{Centre for Mathematics and Its Applications, Australian National University, Canberra, ACT 2601, Australia.}
\email{Haodi.Chen@anu.edu.au}

\subjclass[2010]{52A40, 46E22, 35J60}

\keywords{Blaschke-Santal\'o inequality}

\begin{abstract}
In this paper, we establish a generalised Blaschke-Santal\'o inequality for convex bodies in $\R^{n+1}$.
This inequality gives an upper bound estimate
for the product of dual quermassintegrals of convex body and its polar set.
Our argument is based on induction on dimensions.
\end{abstract}

\maketitle

%%%%%%%%%%%%%%%%%%%%%%%%%%%%%%%%%%%%%%%%%%%%%%%%%%%%%%%%
\baselineskip18pt
\parskip3pt

\section{Introduction}

The theory of mixed volumes is at the core of the study of
geometric invariants and geometric measures associated with convex bodies.
It arises from the combination of the two fundamental concepts of Minkowski addition and volume,
and forms a central part of the Brunn-Minkowski theory.
In the Euclidean $\R^{n+1}$, the quermassintegrals, $W_0,\ldots,W_n$, are the elementary mixed volumes
which include volume, surface area, and mean width \cite{S14}.
Let $\K$ be the set of convex bodies, namely the compact and convex subsets in $\R^{n+1}$
with non-empty interior.
For $\Om\in\K$,
\beqs
W_{n+1-i}(\Om) = \frac{\omega_{n+1}}{\omega_i} \int_{G(n+1,i)}\vol_{i}(\Om|\xi)d\xi,
\eeqs
where $G(n+1,i)$ is the Grassmann manifold of $i$-dimensional subspaces in $\R^{n+1}$,
$\Omega|\xi$ is the image of the orthogonal projection of $\Omega$ onto $\xi$,
$\vol_i$ is the $i$-dimensional volume
and $\omega_i$ is the $i$-dimensional volume of the $i$-dimensional unit ball.
The integration is with respect to the rotation-invariant probability measure on $G(n+1,i)$ \cite{HLYZ16}.

There is a dual Brunn-Minkowski theory introduced in 1970s in \cite{L75}.
Let $\K_0$ be the set of convex bodies which contain the origin in their interiors.
The dual quermassintegrals, $\widetilde W_0,\ldots,\widetilde W_n$, of $\Om\in\K_0$,
can be defined by \cite{HLYZ16}
\beq\label{S1E1}
\widetilde W_{n+1-i}(\Omega)
=\frac{\omega_{n+1}}{\omega_{i}}\int_{G(n+1,i)}\vol_{i}(\Omega\cap\xi)d\xi.
\eeq
It is of interest to study the upper and lower bounds of 
\beqs
W_i(\Omega)W_i(\Omega^*) \ \ \text{and} \ \ \widetilde W_i(\Omega)\widetilde W_i(\Omega^*),
\eeqs
where $\Omega^*$ is the polar body of $\Omega\in\K_0$ with respect to the origin.

The volume is both the quermassintegral $W_0$ and dual quermassintegral $\widetilde W_0$ :
$\vol(\Om)=W_0(\Om) = \widetilde W_0(\Om), \ \forall \ \Om\in\K_0$.
The classical Blaschke-Santal\`o inequality \eqref{claBS} below gives a sharp upper bound of
$W_0(\Omega)W_0(\Omega^*)$  and $\widetilde W_0(\Omega)\widetilde W_0(\Omega^*)$,
for $\Om\in\K_0^e$, the set of all origin-symmetric convex bodies.
This inequality can be also applied to non-symmetric convex bodies.
For $z\in \text{int}\Om$, the dual body of $\Om$ with respect to $z$ is defined as:
\[
\Omega_z^*=\Big\{y+z\ \big|\ y\in \R^{n+1},y\cdot(x-z)\leq1\text{ for all }x\in\Omega\Big\}.
\footnote{In this paper, we sometimes omit the subscript ``$z$" if the dual body of $\Om\in\K_0$
is taken with respect to $z=0$, the origin in $\R^{n+1}$.}
\]
The classical Blaschke-Santal\`o inequality states \cite{S14}
\begin{equation}
\label{claBS}
\sup_{\Om\in\K}\inf_{z\in\text{int}\Om}\vol(\Om)\vol(\Om_z^*)\le \vol^2(B_1),
\end{equation} 
where $B_1$ is the unit ball in $\mathbb R^{n+1}$.
On the other hand,
the minimum of $W_0(\Om)W_0(\Om^*)$ or equivalently $\widetilde W_0(\Om)\widetilde W_0(\Om^*)$
among $\Om\in\K^e_0$, known as Mahler's conjecture, is still unsolved except for dimension two \cite{M39}.
A good lower bound estimate was proved by Bourgain-Milman \cite{BM87}.

The above type inequalities for general indices were studied by several authors
and some partial results were established.
For example, the sharp lower bound for $W_n(\Om)W_n(\Om^*)$
and a lower bound estimate for $W_1(\Om)W_1(\Om^*)$
were proved respectively by Lutwak \cite{L75} and Chai-Lee \cite{CL99},
while upper and lower bounds for $\widetilde W_i(\Om)\widetilde W_i(\Om^*)$,
$i=1,\ldots,n$, were also obtained in \cite{G91} and  \cite{CL99} respectively.

For $\Om\in\K_0$, one can rewrite \eqref{S1E1} as 
\beqs
\widetilde W_{n+1-i}(\Omega)=\frac{1}{n+1}\int_{\S^n}r^i d\s,
\eeqs
where $r=r_\Om$ is the radial function of $\Om$ with respect to the origin,
and $\s$ denotes the standard measure on the unit sphere $\S^n$ \cite{L75}.
It is natural to extend the definitions of dual quermassintegrals to all real indices
\beqs
\widetilde W_q (\Omega)=\frac{1}{n+1}\int_{\S^n}r^{n+1-q}d\s, \ \text{for} \ q\in \R.
\eeqs 
Note that the radial function depends on the choice of the centre.
Hence the dual quermassintegral of $\Om$ may differ for different centres unless $q = 0$.
Let $r_z=r_{\Om,z}:\mathbb S^n\to \mathbb R$ be the radial function of $\Omega\in\K$
with respect to a centre $z\in\text{int}\Omega$, i.e.,
\[r_z(x)=\sup\big\{\lambda>0|\lambda x+z \in \Omega\big\}.\]
It is well known that
\[ \widetilde W_0(\Om)=\text{Vol}(\Omega)=\frac{1}{n+1}\int_{\mathbb S^n} r_z^{n+1} d\s .\]
If $r_z^*$ is the radial function of $\Omega_z^*$ with respect to $z$, then \eqref{claBS} becomes 
\begin{equation}
\label{CBS}
\sup_{\Omega\in\K}\inf_{z\in\text{int}\Omega}\Big(\int_{\S^n}r_z^{n+1}d\s\Big)\Big(\int_{\S^n}{r_z^*}^{n+1}d\s\Big)\leq(n+1)^2\text{Vol}^2(B_1).
\end{equation}
The infimum in \eqref{CBS} is attained at the so-called Santal\`o point and the equality holds
if and only if $\Omega$ is an ellipsoid.
The inequality was firstly proved by Blaschke \cite{B23} for $n=1,2$ in 1917,
and was proved later on by Santal\`{o} \cite{S48} for all dimensions in 1948. 
It was also studied by different approaches \cite{BK15,L85,LZ97,MP90,P85,R81}.

Our main purpose in this paper is to obtain a generalisation of the Blaschke-Santal\`o inequality
in the dual Brunn-Minkowski theory,
which gives an upper bound estimate for
$$\widetilde W_{n+1-\alpha}(\Omega)\widetilde W_{n+1-\beta}(\Omega^*_z)$$
when $\alpha$ and $\beta$ are properly related
(i.e. $\alpha$ and $\beta$ satisfy \eqref{main} or equivalently \eqref{main*} below).

\begin{theorem} \label{mainThm}
Let $\alpha,\beta\in\R_+$. If
\begin{equation}\label{main}
\frac{n}{\alpha}+\frac{1}{\beta}\geq 1\ \ \text{and}\ \ \frac{1}{\alpha}+\frac{n}{\beta}\geq 1,
\end{equation}
then there is a constant $C(n,\alpha,\beta)$ such that
\beq\label{BSineq}
\sup_{\Omega\in\K}\inf_{z\in\text{int}\Omega}\Big(\int_{\S^n}r_z^\alpha d\s \Big)^\frac{1}{\alpha}
\Big(\int_{\S^n}{r_z^*}^\beta d\s\Big)^\frac{1}{\beta}\leq C(n,\alpha,\beta).
\eeq
\end{theorem}

The classical Blaschke-Santal\`{o} inequality is a fundamental geometric inequality,
which has been widely used in the theory of convex body, functional analysis, and PDEs.
For example, Chou-Wang \cite{CW06} obtained the existence of solutions to the $L_p$ Minkowski problem
for $p\in(-n-1,1)$ by making use of the classical Blaschke-Santal\`o inequality.
Similarly, \eqref{BSineq} is useful in the study of
$L_p$ dual Minkowski problem introduced by Lutwak-Yang-Zhang \cite{LYZ18} most recently.
%This new family of geometric measures includes the $L_p$ surface area measures \cite{Lut93}
%and the dual curvature measures \cite{HLYZ16} as special cases,
%and therefore partially unifies the theory of mixed volumes and dual mixed volumes.
%The  $L_p$ dual Minkowski problem \cite{LYZ18} is equivalent to solving
This problem is equivalent to solving the following Monge-Amp\`ere type equation of the support function $u$,
\begin{equation} \label{MAeq}
\det (\nabla^2u+uI)=u^{p-1}\sqrt{u^2+|\nabla u|^2}^{n+1-q} f \ \text{on} \ \S^n,
\end{equation}
where $\nabla$ is the derivative with respect to an orthonormal frame on $\S^n$.
%In fact $\Om\in\K_0$ with uniformly convex and $C^2$ boundary solves the $L_p$ dual Minkowski problem
%with the prescribed measure $d\mu = f d\s$ if and only if its support function $u=u_\Om$ solves \eqref{MAeq},
%where $u:\S^n\to\R$ is given by
%\[
%u(x)=\sup\{\left<x,y\right>:\ y\in \Omega\}.
%\] 
When $p=0, q\in\R$, \eqref{MAeq} is the dual Minkowski problem
proposed by Huang-Lutwak-Yang-Zhang \cite{HLYZ16},
which has been studied by many authors \cite{BHP18,BLYZZ18,CL17,HZ18, LSW17, Zh17a}.
When $p\in\R, q=n+1$, \eqref{MAeq} is the $L_p$ Minkowski problem introduced by Lutwak \cite{Lut93}
and has received intensively investigated,
see e.g. \cite{BLYZ13,CLZ17a,CLZ17b, CW06, Lut-Ol95} and the references therein.
For general $p,q\in\R$, some partial results for the existence of solutions to \eqref{MAeq}
were obtained in \cite{BF18,HZ18}.
In a subsequent paper \cite{CCL18},
the generalised Blaschke-Santal\`o inequality \eqref{BSineq} will be applied to
show that \eqref{MAeq} admits solutions for a large class of $p$ and $q$ under the symmetric assumption.

As the main motivation for the generalised Blaschke-Santal\`o inequality in this paper
is to study the $L_p$ dual Minkowski problem,
we do not pursue any sharp estimates in \eqref{BSineq}.
Certainly, the study of the optimal bound and the equality case for \eqref{BSineq}
are of significant interest and importance.
However, we would like to point out that the relation of $\alpha$ and $\beta$ in Theorem \ref{mainThm} is sharp.

\begin{remark}\label{counter}
In general, \eqref{BSineq} fails if \eqref{main} is not satisfied.
For this, an example will be presented after the proof of Theorem \ref{mainThm}.
\end{remark}

\begin{remark}\label{exponents}
Given $\alpha>0$, let
\beq\label{a*}
\alpha^*=\left\{
\begin{split}
\frac{\alpha}{\alpha-n}\ \ \ \ &\text{if}\ \ \alpha>n+1,\\
n+1\ \ \ \ &\text{if}\ \ \alpha=n+1,\\
\frac{n\alpha}{\alpha-1}\ \ \ \ &\text{if} \ \ 1<\alpha<n+1,\\
+\infty\ \ \ \ &\text{if}\ \ 0<\alpha\leq1.
\end{split}
\right.
\eeq
It is not hard to verify that \eqref{main} is equivalent to
\beq\label{main*}
\alpha>0, \ \text{and} \ \beta\in(0,\alpha^*], \ \beta\ne +\infty.
\eeq
\end{remark}

Our idea of the proof of Theorem \ref{mainThm} is as follows.
By John's lemma, for $\Om \in \K$, there is a $z_e\in\Om$ and an ellipsoid $E$ centred at the origin
such that
\beqs
E+z_e\subset\Om \subset (n+1)E+z_e.
\eeqs
As we take infimum with respect to $z\in\Om$ in \eqref{BSineq},
it is not hard to see that we only need to prove Theorem \ref{mainThm}
for origin-symmetric convex bodies and their dual bodies with respect to the origin.
Namely it suffices to prove Theorem \ref{mainThm OS} stated below.
Since each $\Om\in\K_0^e$ can be compared with rhombi $D=D(a_1,\cdots,a_{n+1})$
(see the definition in Section 2),
the proof of Theorem \ref{mainThm OS} then reduces to properly estimate
the upper bound of $\int_{\S^n} r^\alpha_{D} d\s$ in terms of the parameters $\{a_i\}_{i=1}^{n+1}$.
For $n=1$, such estimate is immediate.
For high dimensions, such estimate is based on a delicate induction argument.

%\begin{remark}
%In terms of $\widetilde W_i(\Omega)\widetilde W_j(\Omega^*)$,
%our inequality \eqref{BSineq} gives an upper bound estimate for $i,j<n+1$ satisfying 
%\beq
%j\geq\left\{
%\begin{split}
%-\frac{ni}{1-i}\ \ \ \ &\text{if}\ \ i<0,\\
%0\ \ \ \ &\text{if}\ \ i=0,\\
%-\frac{i}{n-i}\ \ \ \ &\text{if}\ \ 0<i<n,\\
%-\infty\ \ \ \ &\text{if}\ \ n\leq i<n+1.
%\end{split}
%\right.
%\eeq
%\end{remark}
%\vskip20pt
% \begin{equation}
% \begin{aligned}
% 0<q\leq1,\ \ \ \ \ \ \ \ \ \ \ \ \ \ \ \ \ \ \ \  &p< 0,\\
% 1<q<n+1, \ \ -\frac{nq}{q-1}< &p<0,\\
% n+1\leq q,\ \ \ \ \ \ \  -\frac{q}{q-n}< &p<0.
% \end{aligned}
% \end{equation}

\bigskip

\noindent{\bf Acknowledgement}
The author would like to thank the referee for his/her suggestion on
the improvement of the presentation of the paper.

%----------------------------------------------------------------------------------

\section{Proof of Theorem \ref{mainThm}}

We first give some notations.
Let $a_i$, $i=1,\ldots,n+1$, be positive numbers.
Denote by $D(a_1,\cdots,a_{n+1})$ the rhombus in $\R^{n+1}$ centred at the origin,
with vertices at $(\pm a_1,0,\cdots,0)$, $(0,\pm a_2,0,\cdots,0)$, $\cdots$, $(0,\cdots,0,\pm a_{n+1})$,
by $R(a_1,\cdots,a_{n+1})$ the rectangle with vertices at $(\pm a_1,\cdots,\pm a_{n+1})$,
and by $E(a_1,\cdots,a_{n+1})$ the ellipsoid
\beqs
\Big\{x\in\mathbb R^{n+1}\big| \sum_{i=1}^{n+1}x_i^2/a_i^2\leq1\Big\}.
\eeqs

Let $0$ be the origin.
One can easily verify that 
\begin{equation}\label{dual}
\begin{aligned}
D_0^*(a_1,\cdots,a_{n+1})&=R(\frac{1}{a_1},\cdots,\frac{1}{a_{n+1}}):=R^{-1}(a_1,\cdots,a_{n+1}),\\
R_0^*(a_1,\cdots,a_{n+1})&=D(\frac{1}{a_1},\cdots,\frac{1}{a_{n+1}}):=D^{-1}(a_1,\cdots,a_{n+1}),\\
E_0^*(a_1,\cdots,a_{n+1})&=E(\frac{1}{a_1},\cdots,\frac{1}{a_{n+1}}):=E^{-1}(a_1,\cdots,a_{n+1}),
\end{aligned}
\end{equation}
and
\beqn\label{covering}
D(a_1,\cdots,a_{n+1})&\subset& E(a_1,\cdots,a_{n+1}) \notag \\
&\subset& R(a_1,\cdots,a_{n+1})\subset(n+1)D(a_1,\cdots,a_{n+1}).
\eeqn
Note that for any $\Om\in\K_0$ we always use the notations: $t\Om= \{ty\in\R^{n+1} \ | \ y\in\Om\}$,
and $\Om+z=\{y+z\in\R^{n+1}: \ y\in\Om\}$.

By John's lemma, for any $\Om\in\K$,
there is an ellipsoid $E=E(a_1,\cdots,a_{n+1})$ and $z_e\in\Om$ such that
\beqs
E\subset\Omega-z_e\subset (n+1)E.
\eeqs
Without loss of generality, we may assume $z_e=0$ by translation.
Let $D=D(a_1,\cdots,a_{n+1})$ be the rhombus. By \eqref{covering}, 
\[
D\subset E\subset (n+1)D.\]
Hence one has 
\[
D\subset\Omega\subset(n+1)^2 D \ \ \text{and}\ \ \frac{1}{(n+1)^2}D^*\subset\Omega^*\subset D^*.
\]
This implies that, for any $\alpha,\beta\in\R_+$,
\beqn\label{final}
\Big(\int_{\S^n} r_\Om^\alpha d\s\Big)^\frac{1}{\alpha}\Big(\int_{\S^n} r_{\Om^*}^\beta d\s\Big)^\frac{1}{\beta}
&\leq&\Big(\int_{\S^n} r_{(n+1)^2D}^\alpha d\s\Big)^\frac{1}{\alpha}
\Big(\int_{\S^n}{r_{D^*}}^\beta d\s\Big)^\frac{1}{\beta} \notag\\
&=&(n+1)^2\Big(\int_{\S^n} r_{D}^\alpha d\s\Big)^\frac{1}{\alpha}
\Big(\int_{\S^n}r_{D^*}^\beta d\s\Big)^\frac{1}{\beta}.
\eeqn
Recall that for \eqref{BSineq} we need to take infimum with respect to $z\in\Om$.
This is obviously controlled by the left hand side of \eqref{final}.
Hence it suffices to prove Theorem \ref{mainThm} for origin-symmetric convex bodies,
and in particular, for rhombi.
Namely it suffices to prove the following theorem.
\begin{theorem}\label{mainThm OS} 
Let $\alpha,\beta\in\R_+$. If $\alpha,\beta$ satisfy \eqref{main} or equivalently \eqref{main*},
then there is a positive constant $C(n,\alpha,\beta)$ such that
\beq\label{BSI OS}
\Big( \int_{\S^n} r_\Om^\alpha d\s \Big)^{\frac1\alpha} \Big( \int_{\S^n} r_{\Om^*}^\beta d\s \Big)^{\frac1\beta} 
\le C(n,\alpha,\beta), \ \ \forall \ \Om\in\K_0^e.
\eeq
\end{theorem}

In order to show \eqref{BSI OS}, by virtue of \eqref{final}, it suffices to prove 
\begin{equation}\label{id}
(\int_{\S^n}r_D^\alpha)^\frac{1}{\alpha}(\int_{\S^n}{r_{D^*}}^\beta)^\frac{1}{\beta}\leq C(n,\alpha,\beta),
\ \ \forall \ D = D(a_1,\ldots,a_{n+1}),
\end{equation}
provided $\alpha,\beta$ satisfy \eqref{main} or equivalently \eqref{main*}.

%Now we prove \eqref{id} by induction on dimensions.
\begin{proof}[Proof of inequality \eqref{id}] 

For $0<\alpha,\beta\leq n+1$, it follows by the H\"older's inequality that
\begin{equation}\label{Holder}
\begin{split}
&\Big(\fint_{\S^n} r_D^\alpha d\s\Big)^\frac{1}{\alpha}\leq \Big(\fint_{\S^n} r_D^{n+1} d\s\Big)^\frac{1}{n+1}, \\
&\Big(\fint_{\S^n} r_{D^*}^\beta d\s\Big)^\frac{1}{\beta}\leq \Big(\fint_{\S^n} r_{D^*}^{n+1}d\s\Big)^\frac{1}{n+1}.
\end{split}
\end{equation}
Hence one obtains \eqref{id} by the classical Blaschke-Santal\'o inequality \eqref{CBS}. 

For $\alpha,\beta>n+1$, one sees that such $\alpha,\beta$ do not satisfy \eqref{main}
and so Blaschke-Santal\'o type inequality cannot hold. See Remark \ref{counter}.

It remains to consider the case $\alpha>n+1$ and $0<\beta\leq n+1$,
for which \eqref{main} is equivalent to 
\[0<\beta\leq\frac{\alpha}{\alpha-n}.\]
By H\"older's inequality, it suffices to prove \eqref{id} when $\beta=\frac{\alpha}{\alpha-n}$.
For $0<\alpha\leq n+1$ and $\beta>n+1$,
we only need to exchange $\alpha$ and $\beta$, and the same argument applies.

For $D=D(a_1,\cdots,a_{n+1})$ and $D^{-1}=D(\frac{1}{a_1},\cdots,\frac{1}{a_{n+1}})$,
we may assume that
\beqs
a_1\geq a_2\geq\cdots\geq a_{n+1}>0.
\eeqs By \eqref{dual} and \eqref{covering} one has
\begin{equation*}
\Big(\int_{\S^n} {r_{D^*}}^\beta d\s\Big)^\frac{1}{\beta}
=\Big(\int_{\S^n} r_{R^{-1}({a_1},\cdots,{a_{n+1}})}^\beta d\s\Big)^\frac{1}{\beta}
\leq C(n)\big(\int_{\S^n} r_{D^{-1}}^\beta d\s\big)^\frac{1}{\beta}.
\end{equation*}
Hence, for \eqref{id}, one needs to evaluate
\beqs
\int_{\S^n} r_D^\alpha d\s \ \ \text{and} \ \ \int_{\S^n} r_{D^{-1}}^\beta d\s.
\eeqs
By the symmetry,
we only need to compute in the first octant.
In Cartesian coordinate, $\partial D^{-1}\cap \{x\in \R^{n+1}|x_i\ge0 \text{ for all }i\}$ can be written as
\[
a_1x_1+a_2x_2+\cdots+a_{n+1}x_{n+1}=1,\]
and in spherical coordinates
\[
ra_1\cos\theta_1+ra_2\sin\theta_1\cos\theta_2+\cdots+ra_n\sin\theta_1\cdots\cos\theta_n+ra_{n+1}\sin\theta_1\cdots\sin\theta_n=1,\] 
where $\theta_1,\cdots,\theta_n\in[0,\frac{\pi}{2}]$. Hence 
% under some assumptions:\\
% By assumption (ii), $\Omega_o^*$ (will be written as $\Omega^*$, if no confusion occurs) is a rectangular-like polyhedron, centred at the origin, and the length of side paralleled to the $i-th$ axis is $\frac{2}{a_i}$, for $i=1,\cdot,n+1$. There exist two 'diamond-like' polyhedron, denoted by $D_1$ and $D_2$, whose vertices are $(\pm\frac{1}{a_1},0,\cdots,0),\cdots,(0,\cdots,0,\pm\frac{1}{a_{n+1}})$ and $(\pm\frac{2}{a_1},0,\cdots,0),\cdots,$ $(0,\cdots,0,\pm\frac{2}{a_{n+1}})$, respectively such that $D_1\subset \Omega^*\subset D_2$. Hence 
% \begin{equation*}
% \int(r_{D_1}^q)^\frac{1}{q}\leq\int(r_{\Omega^*}^q)^\frac{1}{q}\leq\int(r_{D_2}^q)^\frac{1}{q}=C(n)\int(r_{D_1}^q)^\frac{1}{q}.
% \end{equation*}
% It suffices to evaluate $(\int r_\Omega^p)^\frac{1}{p}\cdot(\int r_{D_1}^q)^\frac{1}{q}$ to prove \eqref{id}. In this proof, we use the spherical coordinates for computation and by symmetry, it is sufficient to compute $(\int r_{D_1}^p)^\frac{1}{p}$ in the first hyperoctant. In the first hyperoctant, $M_{D_1}$ can be written as
% \[
% ra_1\cos\theta_1+ra_2\sin\theta_1\cos\theta_2+\cdots+ra_n\sin\theta_1\cdots\sin\theta_{n-1}\cos\theta_n+ra_{n+1}\sin\theta_1\cdots\sin\theta_n=0,
% \]
% where $r$ is the radial function of $D_1$ and $\theta_1,\cdots,\theta_n\in[0,\frac{\pi}{2}]$. Hence
\begin{equation}\label{integral}
\int_{\mathbb S^{n}_+} r_{D^{-1}}^\beta d\s
=\int_{[0,\frac{\pi}{2}]^n}\frac{\sin^{n-1}\theta_1\sin^{n-2}\theta_2\cdots\sin\theta_{n-1}d\theta_1\cdots d\theta_n}
{(a_1\cos\theta_1+a_2\sin\theta_1\cos\theta_2+\cdots+a_{n+1}\sin\theta_1\cdots\sin\theta_n)^\beta},
\end{equation}
where
\[\mathbb S^{n}_+=\mathbb S^n\cap\{x\in \R^{n+1}|x_i\ge0 \text{ for all }i\}.
\]
Once the upper bound for \eqref{integral} is obtained,
we can also control $\int_{\mathbb S^+} r_D^\alpha d\s$ from above
by replacing $\beta,a_{1},\cdots,a_{n+1}$ with $\alpha,\frac{1}{a_{n+1}},\cdots,\frac{1}{a_{1}}$. Denote 
\[S(\beta,a_1,\cdots,a_{n+1})=\text{RHS of \eqref{integral}}.\]
Then for \eqref{id} it suffices to show
\begin{equation}\label{induction}   
S^\frac{1}{\beta}(\beta,a_1,\cdots,a_{n+1})\cdot S^\frac{1}{\alpha}(\alpha,\frac{1}{a_{n+1}},\cdots,\frac{1}{a_1})
\leq C(n,\alpha,\beta).
\end{equation}

We shall first prove \eqref{induction} for $n=1$;
and then derive a good bound for $S(\beta, a_1,\cdots,a_{n+1})$
(as well as $S(\alpha, a_1^{-1},\cdots,a_{n+1}^{-1})$)
for general $n>1$ by an induction argument on dimensions,
from which \eqref{induction} follows.

{ \bf Case: $n=1$.}
Let $a_1=\gamma a_2$, with $\gamma\geq1$.
When $\gamma$ is bounded, say for example $\gamma\leq\sqrt2$, one has 
\beqs
S^\frac{1}{\beta}(\beta,a_1,a_2)\cdot S^\frac{1}{\alpha}(\alpha,\frac{1}{a_2},\frac{1}{a_1})
&=& \Big(\int_0^{\frac{\pi}{2}}\frac{d\theta}{(\gamma a_2\cos\theta+a_2\sin\theta)^\beta}\Big)^\frac{1}{\beta}
\Big(\int_{0}^{\frac{\pi}{2}}\frac{d\theta}{(\frac{1}{a_2}\cos\theta+\frac{1}{\gamma a_2}\sin\theta)^\alpha}
\Big)^\frac{1}{\alpha} \\
%&=&\gamma\Big(\int_{0}^{\frac{\pi}{2}}\frac{d\theta}{(\gamma\cos\theta+\sin\theta)^\beta}\Big)^\frac{1}{\beta}
%\Big(\int_{0}^{\frac{\pi}{2}}\frac{d\theta}{(\gamma\cos\theta+\sin\theta)^\alpha}\Big)^\frac{1}{\alpha}\\
&\leq&C\Big(\int_{0}^{\frac{\pi}{2}}\frac{d\theta}{(\cos\theta+\sin\theta)^\beta})^\frac{1}{\beta}
\Big(\int_{0}^{\frac{\pi}{2}}\frac{d\theta}{(\cos\theta+\sin\theta)^\alpha}\Big)^\frac{1}{\alpha}\\
&\leq& C(\alpha,\beta).
\eeqs

When $\gamma>\sqrt2$, we evaluate $S^\frac{1}{\beta}(\beta,a_1,a_2)$
and $S^\frac{1}{\alpha}(\alpha,\frac{1}{a_2},\frac{1}{a_1})$ separately.
We have
\beqn\label{2DimQ}
S^{1/\beta}(\beta,a_1,a_2)
&=&\frac{1}{a_2}\Big(\int_{0}^{\frac{\pi}{2}}
\frac{d\theta}{(\gamma\cos\theta+\sin\theta)^\beta}\Big)^\frac{1}{\beta} \notag\\
&\leq&\frac{1}{a_2}\Big(\int_0^{\frac{\pi}{2}-\arcsin\frac{1}{\gamma}}\frac{d\theta}{(\gamma\cos\theta)^\beta}
+\int_{\frac{\pi}{2}-\arcsin\frac{1}{\gamma}}^\frac{\pi}{2}\frac{d\theta}{\sin^\beta\theta}\Big)^\frac{1}{\beta} \notag\\
&\leq&\frac{1}{a_2}\Big(\int_0^{\frac{\pi}{2}-\arcsin\frac{1}{\gamma}}\frac{d\theta}{(\gamma\cos\theta)^\beta}
+\int_{\frac{\pi}{2}-\arcsin\frac{1}{\gamma}}^\frac{\pi}{2}\frac{d\theta}{\sin^\beta(\pi/4)}\Big)^\frac{1}{\beta} \notag\\
&\leq&\frac{C(\beta)}{a_2}\Big(\frac{1}{\gamma^\beta}\int_\frac{1}{\gamma}^1\frac{1}{t^\beta} dt
+\gamma^{-1}\Big)^\frac{1}{\beta} \notag \\
&\leq&\begin{cases}
C(\beta)a_2^{-1}\gamma^{-1} &\beta<1;\cr
C(\beta)a_2^{-1}\gamma^{-1}\ln\gamma &\beta=1;\cr
C(\beta)a_2^{-1}\gamma^{-\frac{1}{\beta}} &\beta>1.
\end{cases}
\eeqn
Note that in our case $\beta = \frac{\alpha}{\alpha - n}>1$.
But for the use of the induction argument below,
we have to estimate $S(\beta,a_1,\cdots,a_{n+1})$ for all possible $\beta$.
Similarly,
\beqn\label{2DimP}
S^{1/\alpha}(\alpha,\frac{1}{a_2},\frac{1}{a_1})
&=&\gamma a_2\Big(\int_{0}^{\frac{\pi}{2}}
\frac{d\theta}{(\gamma\cos\theta+\sin\theta)^\alpha}\Big)^\frac{1}{\alpha} \notag \\
&\leq&
\begin{cases}
C(\alpha)a_2&\alpha<1,\cr
C(\alpha)a_2\ln\gamma &\alpha=1\cr
C(\alpha)a_2\gamma^{1-\frac{1}{\alpha}}&\alpha>1.
\end{cases}
\eeqn
Combining \eqref{2DimQ} and \eqref{2DimP}, one concludes when $\alpha> n+1$ (with $n=1$), 
\beqn\label{2Dcase}
S^\frac{1}{\beta}(\beta,a_1,a_2)\cdot S^\frac{1}{\alpha}(\alpha,\frac{1}{a_2},\frac{1}{a_1})
&\leq&
\begin{cases}
C(\alpha,\beta)\gamma^{-\frac{1}{\alpha}}&\beta<1,\cr
C(\alpha,\beta)\gamma^{-\frac{1}{\alpha}}\ln\gamma&\beta=1,\cr
C(\alpha,\beta)\gamma^{1-\frac{1}{\alpha}-\frac{1}{\beta}}&\beta>1.
\end{cases}
\eeqn
Since $\beta=\frac{\alpha}{\alpha-1}$ (for $n=1$), we have $\frac{1}{\alpha}+\frac{1}{\beta}=1 $.
Therefore \eqref{induction} follows from \eqref{2Dcase}.

\vskip10pt

{\bf Case: $n>1$.}
Assume $a_i=\gamma_ia_{i+1}$ with $\gamma_i \ge 1$ for $i\in\{1,2,\cdots,n\}$.
Since 
 \[
 S(\beta,a_1,a_2,\cdots,a_{n+1})\leq S(\beta,a_1,\frac{a_2}{\sqrt2},\cdots,\frac{a_{n+1}}{(\sqrt2)^n})
 \]
and 
 \beqs
 S(\alpha,\frac{1}{a_{n+1}},\frac{1}{a_{n}},\cdots,\frac{1}{a_{1}})
 &=&(\sqrt2)^{-n}S(\alpha,\frac{(\sqrt2)^n}{a_{n+1}},\frac{(\sqrt2)^n}{a_n},\cdots,\frac{(\sqrt2)^n}{a_{1}})\\
 &\le&(\sqrt2)^{-n}S(\alpha,\frac{(\sqrt2^n)}{a_{n+1}},\frac{(\sqrt2)^{n-1}}{a_n},\cdots,\frac{1}{a_1}),
 \eeqs
 we have
 \begin{equation}\label{beta}
 \begin{aligned}
 S^\frac{1}{\beta}(\beta,a_1,\cdots,a_{n+1})&\cdot S^\frac{1}{\alpha}(\alpha,\frac{1}{a_1},\cdots,\frac{1}{a_{n+1}})\\
 %&\leq C(n)S^\frac{1}{\beta}(\beta,a_1,\frac{a_2}{\sqrt2},\cdots,\frac{a_{n+1}}{(\sqrt2)^n})\cdot S^\frac{1}{\alpha}
 %(\alpha,\frac{(\sqrt2^n)}{a_{n+1}},\frac{(\sqrt2)^{n-1}}{a_n},\cdots,\frac{1}{a_1}) \\
 &\le C(n)S^\frac{1}{\beta}(\beta,\bar a_1,\bar a_2,\cdots,\bar a_{n+1})\cdot
 S^\frac{1}{\alpha}(\alpha,\frac{1}{\bar a_{n+1}},\frac{1}{\bar a_n},\cdots,\frac{1}{\bar a_1}),
 \end{aligned}
 \end{equation}
where $\bar a_i := a_i/(\sqrt 2)^{i-1}$.
Note that $\bar a_{i}/ \bar a_{i+1} \ge \sqrt 2 $.
Hence we may directly assume $\gamma_i\geq\sqrt2$ for all $i$ in the sequel.

The key ingredient in our proof is to derive the estimates \eqref{qest} and \eqref{1result} below:
\begin{equation}\label{qest}
S(\beta,a_1,\cdots,a_{n+1})
\leq\begin{cases}
C(n,\beta)a_1^{-\beta}&\beta<1, \cr
C(n,\beta)a_1^{-1}a_2^{-\beta+1}&1<\beta<2, \cr
C(n,\beta)a_1^{-1}a_2^{-1}a_3^{-\beta+2}&2<\beta<3, \cr
\vdots\cr
C(n,\beta)a_1^{-1}a_2^{-1}\cdots a_n^{-\beta+n-1}&n-1<\beta< n, \cr
C(n,\beta)a_1^{-1}a_2^{-1}\cdots a_n^{-1}a_{n+1}^{-\beta+n}&\beta>n.
\end{cases}
\end{equation}
and when $\beta$ is an integer, $1\le \beta\le n$,
\begin{equation}\label{1result}
S(\beta,a_1,\cdots,a_{n+1})\leq C(n,\beta)a_1^{-1}\cdots a_\beta^{-1} \max_{r=\beta,\ldots,n}\{1,\ln\gamma_r\}.
\end{equation}

For $n=1$, \eqref{qest} and \eqref{1result} are exactly \eqref{2DimQ}.
We next prove \eqref{qest} and \eqref{1result} by induction on dimensions.
For this, let us assume that \eqref{qest} and \eqref{1result} hold for all $k=1,\ldots,n-1$.

Let \[M=M(\theta_2,\cdots,\theta_n)=a_2\cos\theta_2+\cdots+a_{n+1}\sin\theta_2\cdots\sin\theta_{n}\] 
and 
 \[\Phi=\Phi(\theta_2,\cdots,\theta_n)=\sin^{n-2}\theta_2\cdots\sin\theta_{n-1}.\]
We have

\beqn\label{nDimQ}
S(\beta,a_1,\cdots,a_{n+1})
&=&\int_{[0,\frac{\pi}{2}]^n}\frac{\sin^{n-1}\theta_1\cdots\sin\theta_{n-1}d\theta_1\cdots d\theta_{n}}
{(a_1\cos\theta_1+\cdots+a_n\sin\theta_1\cdots\cos\theta_n+a_{n+1}\sin\theta_1\cdots\sin\theta_n)^\beta} \notag\\
&=&\int_{[0,\frac{\pi}{2}]^n}\frac{\sin^{n-1}\theta_1\cdot\Phi}{(a_1\cos\theta_1+\sin\theta_1\cdot M)^\beta}
d\theta_1\cdots d\theta_n \notag\\
&\leq&\int_{[0,\frac{\pi}{2}]^{n-1}}\bigg(\int_0^{\frac{\pi}{2}-\arcsin{\frac{1}{\gamma_1}}}
\frac{\sin^{n-1}\theta_1\cdot\Phi}{(a_1\cos\theta_1)^\beta}d\theta_1\bigg)d\theta_2\cdots d\theta_n \notag\\
&&+\int_{[0,\frac{\pi}{2}]^{n-1}}\bigg(\int_{\frac{\pi}{2}-\arcsin{\frac{1}{\gamma_1}}}^\frac{\pi}{2}
\frac{\sin^{n-1}\theta_1\cdot\Phi}{(a_1\cos\theta_1+\sin\theta_1\cdot M)^\beta}d\theta_1\bigg)d\theta_2\cdots d\theta_n \notag\\
&=:&\text{I}+\text{II}.
\eeqn
Similarly as for \eqref{2DimQ}, we obtain
\beqn\label{nDimQ1}
\text{I}&=&
\int_{[0,\frac{\pi}{2}]^{n-1}}\bigg(\int_0^{\frac{\pi}{2}-\arcsin\frac{1}{\gamma_1}}
\frac{\sin^{n-1}\theta_1}{(a_1\cos\theta_1)^\beta}d\theta_1\bigg)\Phi d\theta_2\cdots d\theta_n \notag \\
&\le&\begin{cases}
C(n,\beta)a_1^{-\beta},&\beta<1,\cr
C(n,\beta)a_1^{-\beta}\gamma_1^{\beta-1},&\beta>1.\cr
\end{cases} \notag \\
&=& \begin{cases}
C(n,\beta)a_1^{-\beta},&\beta<1,\cr
C(n,\beta)a_1^{-1}a_2^{1-\beta},&\beta>1.\cr
\end{cases}
\eeqn
%Since $\gamma_1\geq\sqrt2$, we know that $\arcsin\frac{1}{\gamma_1}\in[0,\frac{\pi}{4}]$
%and $\sin\theta_1\in[\frac{\sqrt2}{2},1]$.
On the other hand,
\beqs
\text{II}&\le& C(n) \int_{[0,\frac{\pi}{2}]^{n-1}}
\bigg(\int_{\frac{\pi}{2}-\arcsin\frac{1}{\gamma_1}}^\frac{\pi}{2}
\frac{d\theta_1}{(a_1\cos\theta_1+M)^\beta}\bigg) \Phi d\theta_2\cdots d\theta_n\\
&\le& C(n,\beta) a_1^{-1} \int_{[0,\frac{\pi}{2}]^{n-1}}
\bigg(\int_{M}^{\frac{a_1}{\gamma_1}+M}
\frac{dt}{t^\beta}\bigg) \Phi d\theta_2\cdots d\theta_n \\
&\leq& C(n,\beta) a_1^{-1} \int_{[0,\frac{\pi}{2}]^{n-1}}
\big|(a_{2}+M)^{-\beta+1}-M^{-\beta+1}\big|  \Phi d\theta_2\cdots d\theta_n\\
&\leq& C(n,\beta)a_1^{-1}\Big(\int_{[0,\frac{\pi}{2}]^{n-1}}(a_{2}+M)^{-\beta+1} \Phi d\theta_2\cdots d\theta_n
+\int_{[0,\frac{\pi}{2}]^{n-1}}M^{-\beta+1}  \Phi d\theta_2\cdots d\theta_n\Big)\\
&=& C(n,\beta)a_1^{-1}\Big(\int_{[0,\frac{\pi}{2}]^{n-1}}(a_{2}+M)^{-\beta+1}\Phi d\theta_2
\cdots d\theta_n+S(\beta-1,a_2,\cdots,a_{n+1})\Big).
\eeqs
By $a_2\leq a_2+M\leq (n+1)a_2$, we further deduce
\begin{equation}
\label{nDimQ2}
\text{II}\leq
C(n,\beta) \max\big\{a_1^{-1}a_{2}^{-\beta+1},a_1^{-1}S(\beta-1,a_2,\cdots,a_{n+1})\big\}. 
\end{equation}
Plugging \eqref{nDimQ1} and \eqref{nDimQ2} in \eqref{nDimQ}, we obtain
\begin{equation} \label{Induc}
S(\beta,a_1,\cdots,a_{n+1})
\leq 
\begin{cases}
C(n,\beta)\max\big\{a_1^{-\beta},a_1^{-1}S(\beta-1,a_2,\cdots,a_{n+1})\big\},&\beta<1,\cr
C(n,\beta)\max\big\{a_1^{-1}a_{2}^{-\beta+1},a_1^{-1}S(\beta-1,a_2,\cdots,a_{n+1})\big\}, &\beta>1.
\end{cases}
\end{equation}
It is then not hard to verify \eqref{qest} by our induction assumption.

For $\beta=1$, one has
\beq \label{nDim1}
\text{I} = \int_{[0,\frac{\pi}{2}]^{n-1}}\bigg(\int_0^{\frac{\pi}{2}-\arcsin\frac{1}{\gamma_1}}
\frac{\sin^{n-1}\theta_1}{a_1\cos\theta_1}d\theta_1\bigg) \Phi d\theta_2\cdots d\theta_n 
\le C(n)a_1^{-1}\ln\gamma_1,
\eeq
and
\beqn\label{nDim2}
\text{II}&\le &C(n)\int_{[0,\frac{\pi}{2}]^{n-1}}\bigg(\int_{\frac{\pi}{2}-\arcsin\frac{1}{\gamma_1}}^\frac{\pi}{2}
\frac{d\theta_1}{a_1\cos\theta_1+M}\bigg)\Phi d\theta_2\cdots d\theta_n \notag\\
&\leq& C(n)a_1^{-1}\int_{[0,\frac{\pi}{2}]^{n-1}}\ln\frac{a_2+M}{M}\Phi d\theta_2\cdots d\theta_n \notag\\
&\leq& C(n)a_1^{-1}\int_{[0,\frac{\pi}{2}]^{n-1}}\frac{a_2+M}{M} \Phi d\theta_2\cdots d\theta_n \notag\\
&\leq& C(n)a_1^{-1}+C(n)a_1^{-1}a_2\cdot S(1,a_2,\cdots,a_{n+1}) \notag \\
&\leq& C(n)\max\big\{a_1^{-1},a_1^{-1}a_2\cdot S(1,a_2,\cdots,a_{n+1})\big\}.
\eeqn
Combining \eqref{nDim1} and \eqref{nDim2}, we obtain
\beq\label{1Induc}
S(1,a_1,\cdots,a_{n+1})\leq C(n)\max\{a_1^{-1}\ln\gamma_1, a_1^{-1}a_2\cdot S(1,a_2,\cdots,a_{n+1})\}.
\eeq
Hence \eqref{1result} (with $\beta=1$) follows by our induction assumption.

When $\beta$ is an integer and $2\le \beta\le n$,
we have by \eqref{Induc}
\beqs
S(\beta,a_1,\cdots,a_{n+1})
&\le& C(n,\beta)  \max\{a_1^{-1} a_2^{-\beta+1} ,a_1^{-1} S(\beta-1,a_2,\cdots,a_{n+1})\} \\
&\le&C(n,\beta)  \max\{a_1^{-1}a_2^{-1}\cdots a_\beta^{-1} ,a_1^{-1} S(\beta-1,a_2,\cdots,a_{n+1})\}.
\eeqs
Thus \eqref{1result} follows by induction.

Replacing $\beta,a_{1},\cdots,a_{n+1}$ by $\alpha,\frac{1}{a_{n+1}},\cdots,\frac{1}{a_{1}}$ in
\eqref{qest} and \eqref{1result},
we obtain
\begin{equation}\label{pest}
S(\alpha,\frac{1}{a_{n+1}},\cdots,\frac{1}{a_1})\leq
\begin{cases}
C(n,\alpha)a_{n+1}^{\alpha}&\alpha<1, \cr
C(n,\alpha)a_{n+1}a_{n}^{\alpha-1}&1<\alpha< 2, \cr
C(n,\alpha)a_{n+1}a_{n}a_{n-1}^{\alpha-2}&2<\alpha<3, \cr
\vdots\cr
C(n,\alpha)a_{n+1}a_n\cdots a_2^{\alpha-n+1}&n-1<\alpha< n,\cr
C(n,\alpha)a_{n+1}a_n\cdots a_2a_1^{\alpha-n}&\alpha>n.
\end{cases}
\end{equation}
%and
%\beq\label{1result*}
%S(\alpha,\frac{1}{a_{n+1}},\cdots,\frac{1}{a_1})
%\leq C(n,\beta)a_{n+1}\cdots a_{n-\alpha} \max_{r=\alpha,\ldots,n}\{1,\ln\gamma_r\}.
%\eeq

By virtue of \eqref{qest}, \eqref{1result} and \eqref{pest},
we are able to show \eqref{induction} as follows.
If $\beta=\frac{\alpha}{\alpha-n}$ is not an integer,
we pick an integer $1\le m\le n$ such that $m<\beta<m+1$.
Note that $\alpha>n+1$.
It follows from \eqref{qest} and \eqref{pest} that 
\beqn\label{result}
&&S^\frac{1}{\beta}(\beta,a_1,\cdots,a_{n+1})
\cdot S^\frac{1}{\alpha}(\alpha,\frac{1}{a_{n+1}},\cdots,\frac{1}{a_1})\notag\\
&& \leq C(n,\alpha,\beta)a_1^{-\frac{1}{\beta}}a_2^{-\frac{1}{\beta}}
\cdots a_m^{-\frac{1}{\beta}}a_{m+1}^{-1+\frac{m}{\beta}}
\cdot a_1^{1-\frac{n}{\alpha}}a_2^{\frac{1}{\alpha}}\cdots a_{n+1}^{\frac{1}{\alpha}} \notag \\
&& = C(n,\alpha,\beta) \big(\frac{a_1}{a_2}\big)^{1-\frac{1}{\beta}-\frac{n}{\alpha}}
\big(\frac{a_2}{a_3}\big)^{1-\frac{2}{\beta}-\frac{n-1}{\alpha}}
\cdots\big(\frac{a_m}{a_{m+1}}\big)^{1-\frac{m}{\beta}-\frac{n+1-m}{\alpha}} \\
&&\times\big(\frac{a_{m+1}}{a_{m+2}}\big)^{-\frac{1}{\alpha}}
\cdots\big(\frac{a_{m+1}}{a_{n+1}}\big)^{-\frac{1}{\alpha}}. \notag
\eeqn
Since
\beq\label{very temp}
 1- \frac{k}{\beta}- \frac{n+1-k}{\alpha} \le 0 \ \ \forall \ k=1,\ldots,m,
\eeq
we get \eqref{induction}.

If $\beta=\frac{\alpha}{\alpha-n}$ is an integer, we deduce by \eqref{1result} and \eqref{pest} that
\beqn\label{result2}
&&S^\frac{1}{\beta}(\beta,a_1,\cdots,a_{n+1})
\cdot S^\frac{1}{\alpha}(\alpha,\frac{1}{a_{n+1}},\cdots,\frac{1}{a_1}) \notag\\
&&\leq C(n,\alpha,\beta)a_1^{1-\frac{n}{\alpha}}a_2^\frac{1}{\alpha}\cdots a_{n+1}^\frac{1}{\alpha}
\cdot a_1^{-\frac{1}{\beta}}\cdots a_\beta^{-\frac{1}{\beta}}
\cdot\Big(\max_{r=\beta,\ldots,n}\{1,\ln\gamma_r\}\Big)^\frac{1}{\beta} \notag \\
&&=C(n,\alpha,\beta)\big(\frac{a_1}{a_2}\big)^{1-\frac{n}{\alpha}-\frac{1}{\beta}}
\big(\frac{a_2}{a_3}\big)^{1-\frac{n-1}{\alpha}-\frac{2}{\beta}}
\cdots\big(\frac{a_{\beta-1}}{a_\beta}\big)^{(1-\frac{n-\beta+2}{\alpha}-\frac{\beta-1}{\beta})}  \\
&&\ \ \times\big(\frac{a_\beta}{a_{\beta+1}}\big)^{-\frac{1}{\alpha}}\cdots
\big(\frac{a_\beta}{a_{n+1}}\big)^{-\frac{1}{\alpha}}
\Big(\max_{r=\beta,\ldots,n}\{1,\ln\gamma_r\}\Big)^\frac{1}{\beta} \notag \\
&&\le C(n,\alpha,\beta) \prod_{r=\beta}^n \gamma_r^{-\frac1\alpha}
\times \Big( \max_{r=\beta,\ldots,n}\{1,\ln\gamma_r\}\Big)^{\frac1\beta}. \notag
\eeqn
The last inequality is due to \eqref{very temp}.
Since for fixed $\alpha$ and $\beta$,
$(\ln\gamma_r)^\frac{1}{\beta}/(\gamma_r)^\frac{1}{\alpha}$ are uniformly bounded,
we complete the proof.

\end{proof}

Once inequality \eqref{id} is established,
Theorems \ref{mainThm} \& \ref{mainThm OS} follow immediately.

\vskip 15pt
{\bf{An example for Remark \ref{counter}}}
Suppose $\alpha\geq \beta$.
If not, we simply exchange $\alpha$ and $\beta$.
In this case, $\frac{n}{\alpha}+\frac{1}{\beta}\leq \frac{1}{\alpha}+\frac{n}{\beta}$.
If $\alpha,\beta$ do not satisfy \eqref{main},
one have that \[\frac{n}{\alpha}+\frac{1}{\beta}<1.\]
As in the proof, consider the rhombus $D=D(a_1,a_2,\cdots,a_{n+1})$  in $ \R^{n+1}$.
The polar set of $D$ (with respect to the origin) is the rectangle $R(a_1^{-1},\cdots,a_{n+1}^{-1})$
which contains $D^{-1} = D(a_1^{-1},\cdots,a_{n+1}^{-1})$.
For $z\in\text{int}\Omega$, we have 
\begin{equation*}
\int_{\mathbb S^n}r_D^\alpha d\s
=2^{n+1}\int_{\mathbb S^n_+}r_D^\alpha d\s
\leq\ 2^{n+1}\int_{\mathbb S^n} r_{D,z}^\alpha d\s,
\end{equation*}
and 
\[
\int_{\mathbb S^n} r_{D^{-1}}^\beta d\s=2^{n+1}\int_{\mathbb S^n_+}r_{D^{-1}}^\beta d\s
\leq 2^{n+1}\int_{\mathbb S^n}r_{D_z^*,z}^\beta d\s.
\]
Hence 
\[
 \Big(\int_{\S^n}r_D^{\alpha} d\s\Big)^\frac{1}{\alpha}
 \Big(\int_{\S^n}r_{D^{-1}}^{\beta}\Big)^\frac{1}{\beta}
 \leq C \inf_{z\in\text{int}\Omega}\Big(\int_{\S^n}r_{D,z}^{\alpha}d\s\Big)^\frac{1}{\alpha}
 \Big(\int_{\S^n}r_{D_z^*,z}^{\beta} d\s\Big)^\frac{1}{\beta}.
\]
From \eqref{result} and \eqref{result2} one knows that if
we let $a_2=\cdots=a_{n+1}$ and $\gamma_1=\frac{a_1}{a_2}\to+\infty$,
\begin{equation*}
\Big(\int_{\S^n}r_{D}^{\alpha}\Big)^\frac{1}{\alpha}\Big(\int_{\S^n}r_{D^{-1}}^{\beta}\Big)^\frac{1}{\beta}
\sim O(1)\gamma_1^{1-\frac{n}{\alpha}-\frac{1}{\beta}}\to+\infty\ \ \ \text{as}\ \ \gamma_1\to +\infty.
\end{equation*}
Hence the generalised Blaschke-Santal\'o inequality \eqref{BSineq} fails for such $\alpha$ and $\beta$.

%%%%%%%%%%%%%%%%%%%%%%%%%%%%%%%%%%%%%%%%%%%%%%%%%%%%%%%

\bigskip

\bigskip

%%%%%%%%%%%%%%%%%%%%%%%%%%%%%%%%%%%%%%%%%%%%%%%%%%%%%%%

\end{document}